\documentclass[a4paper,10pt]{report}
\usepackage{amssymb}
\usepackage{latexsym}
\usepackage{natbib}
\usepackage{amsmath}

\def\lemma{\textbf{Lemma}}

\def\corollary{\textbf{Corollary}}

\def\propn{\textbf{Proposition}}
\def\pf{\textbf{Proof }}
\def\thm{\textbf{Theorem}}

\def\cor{\textbf{Corollary}}

\def\H2{\B{H}^{2}}
\def\diamond{\diamondsuit}

\newcommand{\lra}{\longrightarrow}
\newcommand{\nt}{\newtheorem}
\newcommand{\tit}{\textit}
\newcommand{\tbf}{\textbf}

\newcommand{\B}{\Bbb}

\begin{document}

\begin{center}
\tbf{Tightness and computing distances in the curve complex.}

Kenneth J. Shackleton

School of Mathematics, University of Southampton,
Highfield, Southampton, SO17 1BJ, Great Britain.

e-mail: kjs@maths.soton.ac.uk

http://www.maths.soton.ac.uk/$\sim$kjs

[First draft: June 2004; Revised: September 2004]
\end{center}

Abstract: We give explicit bounds on the intersection
  number between any curve on a tight multigeodesic and the two ending
  curves. We use this to construct all tight multigeodesics and so
  conclude that distances are computable. The algorithm applies to all
  surfaces. We recover the finiteness result of
  Masur-Minsky for tight geodesics. The central argument makes no use
  of the geometric limit
  arguments seen in the recent work of Bowditch (2003) and
  Masur-Minsky (2000), and is enough to deduce a computable version of the acylindricity theorem of Bowditch.\\


\tbf{\S$0$. Introduction.} Let $\Sigma$ be a closed connected oriented
surface and let $\Pi \subseteq \Sigma$ be a finite subset. In
[Harv], Harvey associates to the pair $(\Sigma,\Pi)$ a simplicial
complex $\mathcal{C}(\Sigma,\Pi)$ called the curve complex. This
is defined as follows. We shall say that an embedded loop in
$\Sigma-\Pi$ is trivial if it bounds a disc and peripheral if it
bounds a once punctured disc. Let $X=X(\Sigma,\Pi)$ be the set of
all free homotopy classes of non-trivial and non-peripheral
embedded loops in $\Sigma-\Pi$. The elements of $X$ will be
referred to as curves. We take $X$ to be the set of vertices and
deem a family of distinct curves $\{\gamma_{0},\gamma_{1},\ldots,
\gamma_{k}\}$ to span a $k$-simplex if any two curves can be
disjointly realised in $\Sigma-\Pi$. The mapping class group has
cocompact simplicial action on $\mathcal{C}$. This has been
exploited by various authors, see
for example [BF], [Hare] and [Iva].

With the exception of only a few cases, namely $\Sigma$ is a
$2$-sphere and $|\Pi| \leq 4$ and $\Sigma$ is a torus and $|\Pi|
\leq 1$, $X$ is non-empty and the curve complex is connected.
For these non-exceptional cases, it can be verified that the
simplicial dimension $C$ of $\mathcal{C}$ is equal to
$3(genus(\Sigma)-1)+|\Pi|-1$. We see that $(\Sigma,\Pi)$ is in fact
non-exceptional if and only if $C(\Sigma,\Pi)>0$. In each subsequent
section, it is to be assumed that $(\Sigma,\Pi)$ is
such that $C(\Sigma,\Pi)>0$.

When $C(\Sigma,\Pi)>0$ the curve complex can be endowed with a
path-metric by declaring all edge lengths to be equal to $1$. All
that is important here is the $1$-skeleton $\mathcal{G}$
of $\mathcal{C}$, with which $\mathcal{C}$ is quasi-isometric. The induced
metric on $\mathcal{G}$ is known to be unbounded and hyperbolic in the
sense of Gromov [MaMi1], [Bow1]. The boundary of $\mathcal{G}$ is
homeomorphic to the space of minimal geodesic laminations filling $\Sigma-\Pi$,
given any hyperbolic metric on $\Sigma-\Pi$,
endowed with the ``measure forgetting'' topology [Kla], [Ham]. The curve
complex plays a key role in Minsky et al's approach to Thurston's
ending lamination conjecture. A second approach has been proposed by
Rees [R].

All this at first sight suggests that
we may apply the methods of hyperbolic groups and spaces to study
various groups acting on $\mathcal{G}$, in particular the mapping
class group and its subgroups. 
The curve graph, though, is not locally finite or
even fine: As early as the $2$-ball around any vertex of $\mathcal{G}$
these problems are manifest. Tight multigeodesics, introduced in
[MaMi2] and further studied in [Bow2] and [Bow3], address this
problem. Their introduction has been fruitful: Masur and Minsky used
these to study the conjugacy problem in the mapping class group
and Bowditch used these to describe the action of the mapping
class group on the curve complex.

Masur and Minsky [MaMi2] showed that there
are only finitely many tight multigeodesics between any two vertices of the
curve graph and Bowditch [Bow2] improved on this, showing that
there are only uniformly boundedly many curves in any given slice. We
go some way to re-establishing these results, though our bounds depend on the intersection
number of the two ending vertices. Since the arguments given
 here do not rely on passing to geometric limits,
our results can be viewed as addressing the local
finiteness problems as well as offering computability. We see how to
construct geodesics, all tight multigeodesics and compute the distance 
between any two vertices.
These notions of tightness are perhaps stronger than we need.

We introduce the key idea of chords and pulses to measure the interleaving
in the surface of curves lying on a geodesic in $\mathcal{G}$. Pulse
is preserved by the action of the mapping class group.
We establish bounds on
pulse that apply to all geodesics and, combining these
with an appropriate tightness criterion, we establish the finiteness of
tight multigeodesics. These methods are readily applicable to related
complexes.

Lastly, in $\S7$ we use our results to compute stable lengths of all
mapping classes.\\

\tbf{$\S1.$ Tightness in the curve graph and the main results.} Let us remind
ourselves of a few definitions. Associated to any two curves $\alpha$
and $\beta$ is their geometric intersection number
$\iota(\alpha,\beta)$, namely the minimal cardinality of the set $a
\cap b$ among all $a \in \alpha$ and $b \in \beta$. Note
$\iota(\alpha,\alpha)=0$ for all $\alpha$, and $d(\alpha,\beta)
\leq 1$ if and only if $\iota(\alpha,\beta)=0$. For any two curves
$\alpha$ and $\beta$, we have have $d(\alpha,\beta) \leq
\iota(\alpha,\beta) + 1$ (see [Bow1] for a logarithmic bound).

For us, paths in the curve graph shall be sequences of vertices
$\gamma_{0},\gamma_{1},\ldots,\gamma_{n}$ such that $\gamma_{i} \neq
\gamma_{i+1}$ and $\iota(\gamma_{i}, \gamma_{i+1})=0$, that is
$\gamma_{i}$ and $\gamma_{i+1}$ are adjacent, for each $i$. A geodesic
in $\mathcal{G}$ is a distance realising path.

We shall recall the notion of tight multigeodesic due to Bowditch
[Bow2], but that of Masur and Minsky [MaMi2] works equally well
here. Recall that a multicurve is a collection of pairwise distinct
curves of pairwise zero intersection number. Intersection number on
multicurves is defined additively. Recall that a multipath is a
sequence of multicurves $(v_{i})^{n}_{0}$ such that
$d(\gamma_{i},\gamma_{j}) = |i-j|$ for all $\gamma_{i} \in v_{i},
\gamma_{j} \in v_{j}$ and $i < j$. We say that a multipath
$(v_{i})_{0}^{n}$ is tight at $v_{j}$ ($1 \leq j \leq n-1$) if for all
curves $\delta$, whenever $\iota(\delta,v_{j})>0$ we have
$\iota(\delta,v_{j-1})+\iota(\delta,v_{j+1})>0$. We say that
$(v_{i})_{0}^{n}$ is tight if tight at each $v_{j}$ ($1 \leq j \leq
n-1$). A tight multipath $(v_{i})_{0}^{n}$ is a tight multigeodesic if
$d(\gamma_{0},\gamma_{n})=n$ for some (hence any) $\gamma_{0} \in
v_{0}$ and $\gamma_{n} \in v_{n}$. The existence of tight
multigeodesics was established in [MaMi2]. Whether we can always
connect two vertices of $\mathcal{G}$ by a tight geodesic, rather than
having to use multicurves, remains open. 

The main result may be stated as follows.

\nt{A}{\lemma}

\begin{A}There is an explicit increasing function $F : \B{N} \lra
  \B{N}$ such that the following holds. Let $(v_{i})_{0}^{n}$ be any
  multigeodesic tight at $v_{1}$. Then, $\iota(v_{1},v_{n}) \leq
  F(\iota(v_{0},v_{n}))$.
\end{A}

Note that $F(n)$ grows superexponentially with $n$. 
In particular, the loss of the
uniformity of the bounds of [Bow2] appears to be the price of
computability. Even so, these bounds are enough to deduce the visual
connectivity of $\partial \mathcal{G}$ by bi-infinite geodesics and other familiar facts.
Consequences of Lemma 1 include the following.

\nt{1}[A]{\thm}

\begin{1}There exists an explicit algorithm which takes as input
  $\Sigma,\Pi$ and any two curves $\alpha$ and $\beta$ in $\Sigma-\Pi$
  and returns all tight multigeodesics connecting $\alpha$ to $\beta$.
\end{1}

From this we conclude that distances are computable.

\nt{2}[A]{\thm}

\begin{2}There exists an explicit algorithm which takes as input
  $\Sigma,\Pi$ and any two curves $\alpha$ and $\beta$ in $\Sigma-\Pi$
  and returns the distance between $\alpha$ and $\beta$ in
  $\mathcal{G}(\Sigma, \Pi)$.
\end{2}

In his unpublished thesis, J. Leasure [Lea] gives a version of Theorem 3 for closed surfaces of genus at least two. We are grateful to Richard P. Kent IV for alerting us to this.\\


\tbf{Acknowledgements.} The author wishes to thank Brian Bowditch for
introducing him to these questions and for many interesting and helpful
conversations. The author also wishes to thank the EPSRC for
generously supporting his research. Last but not least, he wishes to thank the Max
Planck Insitut f\"ur Mathematik, Bonn, for its hospitality whilst
writing up.\\

\tbf{$\S2$. An overview of the proof to Lemma 1 and the first few
  cases.} 
Neighbouring multicurves on a tight multigeodesic tend to drag
one another round the surface and shield each other from other
curves. Consider any multigeodesic $(v_{i})_{0}^{n}$ and any simple
realisation $c_{i}$ for $v_{i}$, each $i$, such that $c_{i} \cap
c_{i+1} = \emptyset$ for each $i$, $|c_{i} \cap c_{n}| =
\iota(v_{i},v_{n})$ each $i \leq  n-2$ and $c_{i} \cap c_{j} \cap c_{n} =
\emptyset$ for each $i < j \leq n-2$. Suppose that two components $J_{1}$
  and $J_{2}$ of $c_{n} - c_{0}$ are connected by three subarcs of
  $c_{1}$, denoted $g_{1},g_{2}$ and $g_{3}$, that are otherwise
  disjoint from $c_{n}$ and are homotopic relative to
  $c_{n}-c_{0}$. Tightness at $v_{1}$ implies that the ends of at
  least two of these subarcs, say $g_{1}$ and $g_{2}$, are separated
  by a point from $c_{2} \cap c_{n}$. Now $c_{1}$ and $c_{2}$ are
  disjoint so we conclude that there must be a subarc $h$ of $c_{2}$
  connecting $J_{1}$ and $J_{2}$ and sandwiched between $g_{1}$ and
  $g_{2}$. This subarc of $c_{2}$ is shielded from $c_{0}$ by
  $c_{1}$. Indeed, if $g_{1}$ and $g_{2}$ continue from $J_{2}$ and
  both return to $J_{1}$ while remaining homotopic relative to
  $c_{n}-c_{0}$ then $h$ continues to be trapped in between $g_{1}$
  and $g_{2}$ all the way back to $J_{1}$. If the ends of $h$ are not
  separated on $J_{1}$ by a point from $c_{3} \cap c_{n}$ then we may
  surger $h$ along $c_{n}$ to find a new simple loop $c_{2}'$ disjoint
  from both $c_{0}$ and $c_{3}$. In particular, when $\Pi$ is empty
  $c_{2}'$ represents a curve, denoted $\gamma_{2}'$, and we have
  succeeded in finding a new multipath $v_{0},\gamma_{2}',v_{3}$,
  contradicting $d(v_{0},v_{3})=3$. We conclude that the ends of $h$
  must be separated by $c_{3}$.

This analysis continues along $(v_{i})_{0}^{n}$ to higher
indices. Suppose that a long subarc of $c_{n-2}$ is fellow travelled
by two long subarcs of $c_{1}$, one on either side, and kept apart
from $c_{1}$ by long subarcs from $c_{2},c_{3},\ldots,c_{n-3}$ in
turn. Then we may surger $c_{n-2}$ along $c_{n}-c_{0}$ and arrive at a new
curve $\gamma_{n-2}'$ having zero intersection with
$v_{0}$. Furthermore, $\iota(v_{n-1},v_{n})=0$ and so this time we are
guarranteed $\iota(\gamma_{n-2}',v_{n-1}) = 0$. We have succeeded in
finding a multipath connecting $v_{0}$ and $v_{n-1}$ of length less
than $n-1$ and hence we have a contradiction.

We shall see that if we start with a multigeodesic $(v_{i})_{0}^{n}$
which is tight at $v_{1}$ and is such that $\iota(v_{1},v_{n})$ is
large relative to $\iota(v_{0},v_{n})$ (we shall quantify this in
terms of $F$ in $\S5$) then this is exactly the situation we find
ourselves in. Furthermore, the argument only requires tightness at
$v_{1}$.

We now deal with the cases $n=2$ and $n=3$ separately, and here after
assume $n \geq 4$.

\nt{3}[A]{\propn}

\begin{3}For each multigeodesic $v_{0},v_{1},v_{2}$ we have
  $\iota(v_{1},v_{2})=0$. For each multigeodesic
  $v_{0},v_{1},v_{2},v_{3}$ we have $\iota(v_{1},v_{3}) \leq
  2\iota(v_{0},v_{3})$.
\end{3}

\pf We note that for any multigeodesic $u_{0},u_{1},u_{2}$ tight at
$u_{1}$ and any multicurve $z$, we have $\iota(z,u_{1}) \leq
2(\iota(z, u_{0}) + \iota(z, u_{2}))$. In particular, when $z=v_{3}$
we have $\iota(v_{1},v_{3}) \leq 2(\iota(v_{0},v_{3}) +
\iota(v_{2},v_{3})) = 2(\iota(v_{0},v_{3}) + 0) =
2\iota(v_{0},v_{3})$. $\diamond$\\

The same argument fails for $n \geq 4$ since $v_{2}$ and $v_{n}$ are
no longer adjacent.\\

\tbf{$\S3$. The idea of pulse.} We introduce a measure for the
interleaving in the surface $\Sigma-\Pi$ of curves lying on a tight
multigeodesic in $\mathcal{G}(\Sigma,\Pi)$. The same ideas can be
applied to geodesics and multigeodesics tight at a given vertex.

For any given positive integer $n$, let $F_{n}$ denote the free
monoid of rank $n$ generated by the set
$\{e_{1},e_{2},\ldots,e_{n}\}$. Let $R$ denote the relation set
$\{e_{i}e_{j}=e_{j}e_{i} : i < j-1\}$ and let $N$ denote the
congruence on $F_{n}$ generated by $R$. We form the quotient $F_{n}/N$
and refer to the elements of this monoid as \tit{chords}. For any
word $w \in F_{n}$, denote by $w(i)$ the $i$th letter appearing in
$w$. For each $i, j$ define $|e_{i} - e_{j}| = |i - j|$.

Chords naturally arise in the context of paths and multipaths in the
curve graph. Consider a path or multipath $(v_{i})_{0}^{n}$ and choose
simple representatives $c_{i}$ for $v_{i}$ once more, so that $c_{i}
\cap c_{i+1} = \emptyset$ each $i$, $|c_{i} \cap c_{n}| =
\iota(v_{i},v_{n})$ for each $i \leq n-2$ and $c_{i} \cap c_{j} \cap c_{n}
= \emptyset$ for each $i < j \leq n-2$. Let $J$ be any component of
$c_{n}-c_{0}$. Orientate $J$ and use this orientation to ennumerate
the points of $J \cap \bigcup_{0}^{n-2} c_{i}$. This ennumeration
spells out an element $w$ of $F_{n}$ by identifying a point from
$c_{i}$ with the $i$th generator $e_{i}$ of $F_{n}$, each
$i$. Tightness at $v_{i}$ implies that $w$ cannot be of the form
$w=w_{1}e_{i}^{3}w_{2}$, for some $w_{1}, w_{2} \in F_{n}$ and each
$i$. We will consider various subsets of $J \cap
\bigcup_{0}^{n-2}c_{i}$, and ennumerate their elements with the orientation on
$J$ to determine $m$-pulse.

Now suppose that $|w(i) - w(i+1)| >1$. Then $\gamma_{w(i)}$ and
$\gamma_{w(i+1)}$ have non-zero intersection number and so we may
homotop both $c_{w(i)}$ and $c_{w(i+1)}$ near $J$ so as to transpose
the two points of intersection. If we re-ennumerate, we arrive at a
second word $w' \in F_{n}$ with $\overline{w} = \overline{w'}$. In
this way, paths may be viewed as defining chords and tight
multigeodesics pinched chords. Note also that each word in a chord
induced by a path or multipath can be induced by the same path or
multipath, just by considering transpositions.

Let us set about defining the $m$-pulse of a given word and then for a
given chord, for each $2 \leq m \leq n$. For each word $w \in F_{n}$
we consider subwords $u$ satisfying the following
three conditions. Firstly, both the initial and the final letters in
$u$ are equal to $e_{1}$. Secondly, for each $i$ we have
$|u(i)-u(i+1)|\leq 1$. Thirdly, between any two successive $e_{1}$'s
in $u$ there is exactly one $e_{m}$. We define the $m$-pulse of a such
a subword $u$ to be equal to the number of times $e_{m}$ appears in
$u$. We define the $m$-pulse of $w$ to be the maximal $m$-pulse
arising among all such subwords $u$ of $w$ and denote it by
$p_{m}(w)$. Even when $u$ satisfies these criteria and is maximal with
respect to inclusion among all such subwords, it need not realise the
$m$-pulse of $w$.

\nt{4}[A]{\lemma}

\begin{4}Suppose that $v, w \in F_{n}$ represent the same chord, that
  is $\overline{v} = \overline{w}$. Then, there is a natural one-to-one
  correspondence between subwords of $v$ and subwords of $w$. In
  particular, this restricts to a correspondence between subwords of
  $v$ satisfying our three criteria and subwords of $w$ satisfying our
  three criteria and preserves their $m$-pulse, each $2 \leq m \leq n$.
\end{4}

\pf Any two elements of a chord are related by a finite sequence of
transpositions. The result follows by an induction on the length of
such sequences. Note that each transposition fixes every subword
satisfying our three criteria and so preserves $m$-pulse, each $m$. $\diamond$\\

For $2 \leq m \leq n$, we define the \tit{$m$-pulse} of a given chord to be
equal to the $m$-pulse of one (hence any) representative word, and
denote this by $p_{m}(\overline{w})$. We have just
seen that this is well-defined.

Let us complete this section with a few examples and remarks. Chords
may be represented and are determined by words from $\{e_{1},
e_{2},\ldots, e_{n}\}$. For instance,
$e_{1}e_{2}e_{3}e_{4}e_{2}e_{3}e_{4}$ and
$e_{1}e_{2}e_{3}e_{2}e_{4}e_{3}e_{4}$ represent the same chord since
the first $e_{4}$ and the second $e_{2}$ may be transposed. The words
$e_{1}e_{2}e_{1}e_{2}e_{1}$ and
$e_{1}e_{1}\ldots e_{1}e_{2}e_{1}e_{1}\ldots e_{1}e_{2}e_{1}$ 
represent different chords although their
$2$-pulses both equal $2$. The words
$e_{1}e_{1}e_{1}e_{2}e_{1}e_{1}e_{1}e_{2}e_{1}$ and
$e_{1}e_{2}e_{2}e_{2}e_{1}e_{2}e_{2}e_{2}e_{1}$ define different
chords although both their $2$-pulses and their lengths are
equal. When chords are induced by a given geodesic, the tightness
property prevents repetition. If we bound the $m$-pulse on a chord,
each $m$, then we bound its length.

Lastly, note that pulse is symmetric and almost additive: For each
$m$, the $m$-pulse of the concatenation of two words is either the sum
of each $m$-pulse or one more than this sum. The $1$-pulse of a word
or chord should always be regarded as zero and, for each $2 \leq m
\leq n-1$, the $m+1$-pulse of a chord is at most the
$m$-pulse. Summing the $2$-pulses over each component of $c_{n}-c_{0}$
closely approximates $\iota(v_{1},v_{n})$.

\nt{5}[A]{\lemma}

\begin{5}Suppose that $(v_{i})_{0}^{n}$ is a multigeodesic tight at
  $v_{1}$. Then $\sum p_{2}(\overline{w}) \leq \iota(v_{1},v_{n}) \leq 2\sum
  p_{2}(\overline{w}) + \iota(v_{0},v_{n})$, where the summations are taken over the components
  of $c_{n}-c_{0}$.
\end{5}

\tbf{$\S4$. Train tracks relative to the ends of a geodesic in the
  curve graph.} Let us fix a smooth structure on $\Sigma$. A train
  track $\tau$ in $\Sigma-\Pi$ is a smooth branched $1$-submanifold
  such that the Euler characteristic of the smooth double of each
  component of $\Sigma-(\Pi \cup \tau)$ is negative. This rules out
  discs, once-punctured discs and discs with one or two boundary
  singularities as complementary regions. Train tracks were introduced
  by Thurston to study geodesic laminations.

It is standard to refer to the branch points of a train track as
switches and the edges between switches as branches. We say that
$\tau$ is generic if each switch has valence three. By sliding
branches along branches, if need be, we can take a train track and
return a generic train track. This is a convenient option since it
greatly simplifies our counting arguments. A train subpath $p : I \lra
\tau$ is a continuous map on a closed interval $I \subseteq \B{R}$
such that $p(n)$ is a switch for each $n \in \B{Z} \cap I, p^{-1}(v)
\in \B{Z}$ for each switch $v$ and $\partial I \subseteq \B{Z} \cup
\{\pm \infty\}$.

A smooth simple closed loop $c$ is carried by $\tau$ if there exists a
smooth map $\phi : \Sigma-\Pi \lra \Sigma-\Pi$ homotopic to the
identity map on $\Sigma-\Pi$ such that $\phi|_{c}$ is an immersion and
$\phi(c) \subseteq \tau$. We refer to $\phi$ as a carrying map or
supporting map. If $\tau$ carries $c$ then we have a measure on the
branch set of $\tau$ by counting the number of times the subpath
$\phi(c)$ traverses any given branch. This measure satisfies a switch
condition: At each switch, the total inward measure is equal to the
total outward measure.

We recall a useful combinatorial lemma relating the number of switches
and the number of branches of a train track to the Euler
characteristic of
$\Sigma-\Pi$. This is Corollary 1.1.3 from [PenH].

\nt{6}[A]{\lemma}

\begin{6}Let $\tau$ be any train track in $\Sigma-\Pi$, let $s$ denote
  the number of switches and $e$ the number of branches. Then:

i). $s \leq -6\chi(\Sigma-\Pi) - 2|\Pi|$;

ii). $e \leq -9\chi(\Sigma-\Pi) - 3|\Pi|$.
\end{6}

Let $(v_{i})_{0}^{n}$ be any multipath in $\mathcal{G}(\Sigma,\Pi)$
and choose smooth and simple realisations $c_{i}$ for $v_{i}$, each
$i$, such that $c_{i} \cap c_{i+1} = \emptyset$ for each $i$, $c_{i} \cap
c_{n} = \iota(v_{i}, v_{n})$ for each $i \leq n-2$ and $c_{i} \cap
c_{j} \cap c_{n} = \emptyset$ for each $i < j \leq n-2$. We construct a
train track $\tau$ which will carry all of $c_{1}$ and all those
subarcs of each $c_{i}$ ($2 \leq i \leq n-2$) which end on $c_{n}$ and
which are trapped between subarcs of $c_{1}$ over large distances.

There exists a smooth surjection $\phi : \Sigma-\Pi \lra \Sigma-\Pi$
homotopic to the identity map such that the restriction of $\phi$ to
$c_{1}$ is an immersion onto a smooth branched $1$-submanifold $\tau$
of $\Sigma-\Pi$ with the characterising property that any two
components of $c_{1} - c_{n}$ homotopic relative to $c_{n}-c_{0}$ are
carried into the same edge of $\tau$ and $\tau$ is to be disjoint from
$c_{0}$. Each branch point necessarily belongs to one component of
$c_{n}-c_{0}$ and each component of $c_{n}-c_{0}$ contains at most one
branch point. We now check that $\tau$ defines a train track, with each
branch point viewed as a switch and each edge thought of as a branch,
and that $\tau$ is unique up to isotopy.

\nt{7}[A]{\lemma}

\begin{7}$\tau$ is a train track.
\end{7}

\pf Note that no region complementary to $\tau$ can be diffeomorphic
to a disc with smooth boundary or a monogon (disc with one outward
pointing singularity) by the minimality of $|c_{1} \cap c_{n}|$.

Suppose for contradiction that $E$ is a bigon component, that is a disc with two outward
pointing singularities, of $\Sigma-(\Pi \cup \tau)$. The twos subarcs
of $\partial E$ connecting the two singularities of $\partial E$ are
homotopic to one another relative to $c_{n}-c_{1}$. Hence $E$ must
intersect $c_{0}$, for otherwise these two subarcs of $\partial E$ would
have been collapsed into a single branch of $\tau$. Since $\tau$ and
$c_{0}$ are disjoint, so $\partial E$ and $c_{0}$ are disjoint. Hence
$E$ contains a component of $c_{0}$ which is therefore homotopically
trivial. This is absurd, and we conclude that $\tau$ is a train
track. $\diamond$

\nt{8}[A]{\lemma}

\begin{8}Suppose that $c_{0}^{i},c_{1}^{i},\ldots,c_{n}^{i}$
  ($i=1,2$) are two such realisations for
  $v_{0},v_{1},\ldots,v_{n}$ and that $\tau_{1}$ and
  $\tau_{2}$ are the resulting train tracks, respectively. Then
  $\tau_{1}$ and $\tau_{2}$ are isotopic.
\end{8}

\pf This follows since $c_{0}^{1} \cup c_{1}^{1} \cup c_{n}^{1}$ and
$c_{0}^{2} \cup c_{1}^{2} \cup c_{n}^{2}$ are isotopic. $\diamond$\\

It is worth pointing out that the same construction for $c_{i}$
($i=2,3,\ldots,n-2$) will not necessarily yield a train track but
instead a bigon train track, that is we allow the complementary
regions to be bigons. Each complementary bigon will contain at least
one point from $c_{0} \cap c_{n}$ and there would be at most
$\iota(v_{0},v_{n})$ bigons.

To each switch $v$ of $\tau$ we can associate the set $\phi^{-1}(v)
\cap \bigcup_{0}^{n-2} c_{i}$, which we henceforth denote by
$D(v)$. Let us assume that $c_{0},c_{1},\ldots,c_{n}$ are such that
$|D(v)|$ is minimal for each switch $v$. Orientate $c_{n}$ and use
this orientation to ennumerate the points of $D(v)$. This gives us a
word $w$ in $F_{n-2}$. Thus, for each integer $2 \leq m \leq n-2$, we
may associate to the switch $v$ the $m$-pulse of the chord
$\overline{w}$.

We may use pulse on switches to define measures on the branch set of
$\tau$. Suppose that $v_{1}$ and $v_{2}$ are adjacent switches of
$\tau$ connected by a branch $b$. In what follows, the topological closure
of $b$ in $\Sigma-\Pi$ is denoted $cl(b)$.

\nt{9}[A]{\lemma}

\begin{9}For each integer $2 \leq m \leq n-2$, the $m$-pulse of
  $\phi^{-1}(cl(b)) \cap D(v_{1})$ is equal to the $m$-pulse of
  $\phi^{-1}(cl(b)) \cap D(v_{2})$.
\end{9}

\pf For each $x \in \phi^{-1}(cl(b)) \cap D(v_{1})$, define $q(x)$ to
be the end on $\phi^{-1}(cl(b)) \cap D(v_{2})$ of the subarc of the
$c_{i}$ containing $x$, beginning at $x$ and collapsing to $b$. We
denote this arc by $[x,q(x)]$. The arc $[x,q(x)]$ crosses a second arc
$[y,q(y)]$ only if the corresponding multicurves $c_{i}$ and $c_{j}$
satisfy $|i-j| \geq 2$. $\diamond$\\

We define the $m$-pulse of the branch $b$ to be equal to the $m$-pulse
associated to $\phi^{-1}(cl(b)) \cap D(v_{1})$, or indeed the
$m$-pulse associated to $\phi^{-1}(cl(b)) \cap D(v_{2})$, and denote
it by $p_{m}(b)$. This measure on the branch set of $\tau$ satisfies a
coarse switch condition. For any switch $v$ we may choose one of two
directions at $v$ and partition the branches incident on $v$ as
outgoing and incoming. Denote the corresponding branches by
$b_{1},b_{2},\ldots,b_{s}$ and $b_{s+1},b_{s+2},\ldots,b_{s+t}$,
respectively.

\nt{10}[A]{\lemma}

\begin{10}``Coarse switch condition.'' For each integer $2 \leq m \leq
  n-2$, we have:

i). $|\sum_{1}^{s}p_{m}(b_{i}) - \sum_{s+1}^{s+t}p_{m}(b_{j})| \leq
s+t-2$;

ii). $p_{m}(v) - s + 1 \leq \sum_{1}^{s}p_{m}(b_{i}) \leq p_{m}(v)$;

ii'). $p_{m}(v) - t + 1 \leq \sum_{s+1}^{s+t}p_{m}(b_{j}) \leq
p_{m}(v).$
\end{10}

In particular, if $\tau$ is generic we have the following.

\nt{11}[A]{\cor}

\begin{11}Suppose that $\tau$ is generic and that $b_{1}$ and $b_{2}$
  are outgoing. For each $m$ we have:

i). $0 \leq p_{m}(b_{3}) - p_{m}(b_{1}) - p_{m}(b_{2}) \leq 1$;

ii). $p_{m}(v) - 1 \leq p_{m}(b_{1}) + p_{m}(b_{2}) \leq p_{m}(v)$;

iii). $p_{m}(b_{3}) = p_{m}(v)$.
\end{11}

We shall only require Corollary 12i). and we prove it directly.\\

\pf Suppose that $v$ is a trivalent switch with outgoing branches
$b_{1}$ and $b_{2}$ and incoming branch $b_{3}$. Each component of
$c_{1} \cap \phi^{-1}(b_{3})$ either goes on to be supported by
$b_{1}$ or by $b_{2}$. All those components of $c_{m} \cap
\phi^{-1}(b_{3})$ that are trapped between components of
$\phi^{-1}(b_{1} \cup \{v\} \cup b_{i})$ go on to be supported by
$b_{i}$ ($i = 2,3$). However, where these components of $c_{1}$
diverge subarcs of $c_{m}$ may escape. This reduces the total outgoing
$m$-pulse by at most one, if at all. That is, we have part i). of Corollary 11. $\diamond$\\

Each arc $g$ supported by $\tau$ and whose
ends each lie on components of $b-a$ containing a switch of $\tau$ defines a
train subpath $g_{\phi}$ of $\tau$. For a train path $q$ in $\tau$ we
define the \tit{support} $Supp(q)$ of $q$ to be the set of all subarcs
$g$ of each $c_{i}$ ($i=1,2,\ldots,n-2$) with $q = g_{\phi}$. Suppose
that $q : I \lra \tau$ is a train path with $0 \in I \subseteq
[0,\infty)$. We call the branch of $\tau$ containing $q(i + 1/2)$ the
$i$th branch of $q$. For each integer $m \geq 1$ and each $i \geq 0$
we define the $m$-pulse of the $i$th branch of $q$ to be equal to the
$m$-pulse of either end of $\phi^{-1}(q([0,1])) \cap \bigcup
Supp(q|_{[0,i+1]})$. We denote this by $p_{m,q}(i)$. Note Lemma 10
tells us that this is well defined. Note also that $p_{m,q}(0)$ is
precisely the $m$-pulse of the $0$th branch traversed by $q$, each
$m$.

\nt{12}[A]{\lemma}

\begin{12}``Trains run out of fuel.'' Let $q : [0, \infty) \lra \tau$
    be any train path in $\tau$. Then:

i). $p_{m,q}(i+1) \leq p_{m,q}(i)$ for all $m$ and all $i$;

ii). $p_{m,q}(i) \lra 0$ as $i \lra \infty$.
\end{12}

\pf Each time $q$ arrives at a switch of $\tau$, subarcs of $c_{1}$
that have so far induced $q$ may diverge and there can be no gain in
pulse. Hence i). holds. If $p_{m,q}(i) \geq 1$ for all $i$ then, since
$\Sigma$ is orientable, we conclude that $c_{1}$ has two freely
homotopic components and this is absurd. Hence
ii). holds. $\diamond$\\

The sum of all $2$-pulses over each branch of $\tau$ resembles a
reduced intersection number for $c_{1}$ and $c_{n}$ relative to
$c_{0}$. Suppose that $\iota(v_{0},v_{n})$ is large. Then most of the
regions complementary to $c_{0} \cup c_{n}$ are squares. Whenever
$c_{1}$ meets an edge of one of these squares, necessarily from
$c_{n}$, it must go on to meet the edge opposite. Consider a sequence
$S_{1},S_{2}, \ldots, S_{k}$ of closed squares whose interiors are
complementary to $c_{0} \cup c_{n}$ and such that $S_{i} \cap S_{i+1}
\subseteq c_{n}$, each $i$. Whenever $c_{1}$ meets one of the outer
edges of $\bigcup_{1}^{k} S_{i}$ then $c_{1}$ remains trapped in $\bigcup_{1}^{k}
S_{i}$ and goes on to meet each edge $S_{i} \cap S_{i+1}$ before
exiting at the other outer edge. However, the components of $c_{1}
\cap \bigcup_{1}^{k} S_{i}$ all collapse into a single branch of $\tau$. Now
suppose that $S_{1}$ is met by $c_{1}$. If any two such components of
$c_{1}$ are separated by a subarc of $c_{2}$ then the $2$-pulse on
this branch is precisely one less than the number of components of
$c_{1} \cap \bigcup_{1}^{k} S_{i}$.\\

\tbf{$\S5.$ Proof of Lemma 1.} Let $(v_{i})_{0}^{n}$ be any
multigeodesic in $\mathcal{G}(\Sigma,\Pi)$ tight at $v_{1}$ with $n
\geq 4$. Let $c_{i}$ be any realisation for $v_{i}$, each $i$, such that
$c_{i} \cap c_{i+1} = \emptyset$ for each $i \leq n-1$, $|c_{i} \cap
c_{n}| = \iota(v_{i},v_{n})$ each $i \leq n-2$ and
$c_{i} \cap c_{j} \cap c_{n} = \emptyset$ each $i < j \leq
n-2$. Recall the construction of the train track $\tau$ relative to
$v_{0}$ and $v_{n}$ and carrying all of $c_{1}$, as described in
$\S4$. We can endow each branch of $\tau$ with a family of measures
and each of these satisfies a coarse switch condition. For any train
subpath $q$ of $\tau$ we defined a time measure associated to $q$ by
considering those subarcs of each $c_{i}$ that induce $q$ via $\phi$.

We define the function $K_{n} : \{1,2,\ldots,n-1\} \lra \B{N}$ by the
recurrence relation $$K_{n}(j) = 2^{-(9\chi(\Sigma-\Pi) + 3|\Pi|)(1 +
  K_{n}(j+1))},$$ for all $j = 2,3,\ldots, n-2$, with the boundary
  condition $K_{n}(n-1) = 1$ and where $\chi(\Sigma-\Pi)$ denotes the
  Euler characteristic of $\Sigma-\Pi$.

\nt{13}[A]{\lemma}

\begin{13}Let $m \geq 2$. Suppose that all the branches of $\tau$ have
  $m+1$-pulse at most $K_{n}(m+1)$ and that at least one branch of
  $\tau$ has $m$-pulse at least $K_{n}(m)$. Then there exists a subarc
  $h$ of $c_{m}$ disjoint from $c_{0}$ and whose ends lie on, but is otherwise
  disjoint from, one component of $c_{n}-(c_{0} \cup c_{m+1})$.
\end{13}

\pf We remark that $K_{n}(j)$ was chosen to guarrantee the existence
of a circuit $q : [0,k] \lra \tau$ such that $p_{m,q}(k-1) \geq
1$. Hence $Supp(q)$ contains at least two subarcs $g_{1}, g_{2}$ of
$c_{1}$ beginning and ending on the same component of $c_{n} - (c_{0}
\cup c_{m+1})$ and homotopic relative to $c_{n} - c_{0}$ that trap a
subarc $h$ of $c_{m}$ whose ends lie on the same component of
$c_{n}-c_{0}$. $\diamond$

\nt{14}[A]{\cor}

\begin{14}The $m$-pulse on each branch of $\tau$ is at most
  $K_{n}(m)$, each $m \geq 2$.
\end{14}

\pf Suppose, for contradiction, that there exists $m$ and a branch of
$\tau$ whose $m$-pulse is at least $K_{n}(m)$. We take $m$ to be
maximal subject to this property. Now by Lemma 14 we deduce that there
exists a subarc $h$ of $c_{m}$ beginning and ending on the same
component of $c_{n}$ and disjoint from $c_{0}$. When $\Pi$ is empty we
know that the union of $h$ and the subinterval of $c_{n}-c_{0}$
connecting its ends defines a curve, denoted $\delta$, and we know
that this curve has zero intersection with both $v_{0}$ and
$v_{m+1}$. We have found a multipath $v_{0},\delta,v_{m+1}$ of
length two. Since $d(v_{0},v_{m+1}) = m+1 \geq 2 + 1 = 3$ we have a
contradiction.

When $\Pi$ is non-empty we only have to be slightly more careful since
$\delta$ may be peripheral. Instead, if we define $K_{n}'$ by the
recurrence relation $$K_{n}'(j) = 2^{-2(9\chi(\Sigma-\Pi)+3|\Pi|)(1 +
  K_{n}'(j+1))}$$ with the same boundary condition $K_{n}'(n-1)=1$,
then we can ask for the second return to the same component of $c_{n}
- (c_{0} \cup c_{m+1})$. By considering the boundary components of a
regular neighbourhood of the union of $h$ and the subinterval of
$c_{n}-c_{0}$ connecting the ends of $h$, we again find a curve
$\delta$ which again has zero intersection with both $v_{0}$ and
$v_{m+1}$. $\diamond$\\

\nt{15}[A]{\corollary}

\begin{15}Suppose that $(v_{i})_{0}^{n}$ is a multigeodesic tight at
  $v_{1}$. Then $\iota(v_{1},v_{n}) \leq
  2(1+K'_{n}(2))\iota(v_{0},v_{n})$ if $\Pi \neq \emptyset$ and
  $\iota(v_{1},v_{n})\leq 2(1+K_{n}(2))\iota(v_{0},v_{n})$ if
  $\Pi=\emptyset$.
\end{15}

\pf Let us consider $\Pi \neq \emptyset$. We have seen that the
$2$-pulse on each branch of $\tau$ is at most $K_{n}'(2)$. Since
$(v_{i})_{0}^{n}$ is tight at $v_{1}$ we have $\iota(v_{1},v_{n}) \leq
2(1+max\{p_{2}(b) : b$ is a branch of $\tau\})\iota(\alpha,\beta) \leq
(1+K_{n}'(2))\iota(\alpha,\beta)$. $\diamond$\\

We conclude the proof of Lemma 1.\\

\pf (of Lemma 1). Since $d(v_{0},v_{n}) \leq \iota(v_{0},v_{n}) + 1$
and $d(v_{0}, v_{n}) = n$ we have $K_{n}'(2) \leq
K_{\iota(\alpha,\beta)+1}$. Hence $F(s)=(1+K_{s+1}'(2))s$, each $s \in
\mathbb{N}$, suffices. $\diamond$\\

\tbf{$\S6.$ The proof of Theorem 2 and Theorem 3.} In this section, we
prove the main implications of Lemma 1: We establish a finiteness
result for tight multigeodesics and we establish the computability of
distances in the curve graph.

\nt{16}[A]{\lemma}

\begin{16}There exists an explicit increasing function $F_{1} : \B{N}
  \lra \B{N}$ such that the following holds. Let $(v_{i})_{0}^{n}$ be
  any tight multigeodesic. Then, $\iota(v_{j},v_{n}) \leq
  F_{1}(\iota(v_{0},v_{n}))$ for all $j$.
\end{16}

\pf This follows by an inductive argument using Lemma 1, noting that
$\iota(v_{j},v_{n}) \leq F^{j}(\iota(v_{0},v_{n}))$ for $1 \leq j \leq
n-2$. Now $F_{1}$ defined by $F_{1}(k)=F^{k}(k)$ will
suffice. $\diamond$\\

By considering $(v_{n-i})_{0}^{n}$ instead, we deduce the following.

\nt{17}[A]{\cor}

\begin{17}There exists an explicit increasing function $F_{2} : \B{N}
  \lra \B{N}$ such that the following holds. Let $(v_{i})_{0}^{n}$ be
  any tight multigeodesic. Then, $\iota(v_{j},v_{n}) \leq
  F_{2}(\iota(v_{0},v_{n}))$ and $\iota(v_{0},v_{j}) \leq
  F_{2}(\iota(v_{0},v_{n}))$ for all $j$.
\end{17}

Let $\alpha$ and $\beta$ be any two vertices of the curve graph. The
set of multipaths of length at most $\iota(\alpha,\beta)+1$
connecting $\alpha$ to $\beta$ for which each multicurve verifies
the bounds in Corollary 18 has size uniformly and explicitly bounded in terms of
$\iota(\alpha,\beta)$, $genus(\Sigma)$ and $|\Pi|$. In particular,
this set contains all the tight multigeodesics connecting $\alpha$ to
$\beta$ and we are left to look for these in a bounded search
space. This concludes the proof of Theorem 2. Since tight
multigeodesics are distance realising multipaths, we deduce Theorem
3.\\

\tbf{$\S7.$ The computability of stable lengths.} Given a metric space
$X$ and an isometry $h:X \lra X$ we define the stable length, $||h||$,
of $h$ to be equal to $lim_{n \lra \infty} d(x,h^{n}x)/n$. See [BGS]
for more details.
It is
easily verified that that $||h||$ does not depend on the choice of
$x$. We say that a mapping class $h$ is pseudo-Anosov if for any two
curves $\alpha$ and $\beta$ we have $\iota(\alpha,h^{n}(\beta)) \lra
\infty$ as $n \lra \infty$ (see [FLP]).
We consider $\mathcal{G}$ endowed with usual path-metric, and
prove the computability of the stable length of any given
pseudo-Anosov mapping class.

Let us first recall a few results. In [Bow1], not only is the
hyperbolicity of the curve complex re-established 
but it is also shown that one can compute
hyperbolicity constants. In [Bow2], it is established that there exists
a positive integer $N=N(genus(\Sigma),|\Pi|)$ such that for each
pseudo-Anosov mapping class $h$, $h^{N}$ has a geodesic axis in $\mathcal{G}$. This
implies the stable lengths of pseudo-Anosov mapping classes are both positive and
uniformly rational. It is not known how to compute $N$. Lastly, in a
$k$-hyperbolic geodesic metric space each geodesic rectangle is
$8k$-narrow (so that any point on any one side of the rectangle 
is within $8k$ of the
union of the other three). We now state the result in full:

\nt{computstab}[A]{\thm}

\begin{computstab}There is an algorithm which takes as input
  $\Sigma,\Pi, N(\Sigma,\Pi)$ and a pseudo-Anosov mapping class $h$ and returns
  $||h||$.
\end{computstab}

\pf Fix a choice of $k$ such that $\mathcal{G}$ is
$k$-hyperbolic. Choose an integer $M \geq 18k$.
Let us suppose that $h$ is the $Nth$ power of a
pseudo-Anosov. Then, $h$ is again pseudo-Anosov and has a geodesic
axis denoted $L$. Choose any curve $\alpha$ and 
construct a geodesic $[\alpha,h^{M}\alpha]$ from $\alpha$ to
$h^{M}\alpha$ in $\mathcal{G}$. A central vertex $\beta$ of
$[\alpha,h^{M}\alpha]$ 
must lie
within $8k$ of $L$. Now construct a geodesic $[\beta,h^{M}\beta]$ from $\beta$ to
$h^{M}\beta$. We have $M||h|| =||h^{M}|| \leq d(\beta,h^{M}\beta) + 16k \leq
Md(\beta,h\beta) + 16k$. Hence $||h|| \leq d(\beta,h\beta) + 16k/M < d(\beta,h\beta) +
1$. As $||h||$ is an integer, so $||h|| \leq d(\beta,h\beta)$. 

Further, $d(\beta,h^{M}\beta) \leq ||h^{M}|| + 16k = M||h|| + 16k$ and so
$d(\beta,h^{M}\beta)/M \leq ||h|| + 16k/M < ||h|| +1$.

Combining the two inequalities, we have $||h|| \leq d(\beta,h^{M}\beta)/M
< ||h|| + 1$. Hence $\lfloor d(\beta,h^{M}\beta)/M \rfloor =
||h||$. $\diamond$\\

Notice that in the above we do not find an axis $L$. It would
interesting to find a way of doing so.\\

\tbf{References.}

[BF] M. Bestvina, K. Fujiwara, \tit{Bounded cohomology of
subgroups of the mapping class groups} : Geometry \& Topology
\tbf{6} (2002) 69 - 89.

[BGS] W. Ballmann, M. Gromov, V. Schr\"oder, \tit{Manifolds of
  non-positive curvature} : Progress in Mathematics \tbf{61},
  Birkh\"auser (1985).

[Bow1] B. H. Bowditch, \tit{Intersection numbers and the hyperbolicity of the curve complex} : Preprint, Southampton (2002).

[Bow2] B. H. Bowditch, \tit{Tight geodesics in the curve complex} : Preprint, Southampton (2003).

[Bow3] B. H. Bowditch, \tit{Length bounds on curves arising from tight geodesics} : Preprint, Southampton (2003).

[FLP] A. Fathi, F. Laudenbach, V. Po\'enaru, \tit{Travaux de Thurston
  sur les surfaces (seconde \'edition)} : Asterisque \tbf{66-67},
  Soci\'et\'e math\'ematique de France (1991).

[Ham] U. Hamenst\"adt, \tit{Train tracks and the Gromov boundary
of the complex of curves} : Preprint (2004).

[Hare] J. L. Harer, \tit{The virtual cohomological dimension of the mapping class groups of an orientable surface} : Inventiones Mathematicae \tbf{84} (1986) 157-176.

[Harv] W. J. Harvey, \tit{Boundary structure of the modular group} : in ``Riemann surfaces and related topics: Proceedings of the 1978 Stony Brook Conference'' (ed. I. Kra, B. Maskit), Annals of Mathematical Studies No. 97, Princeton University Press (1981) 245-251.

[Iva] N. V. Ivanov, \tit{Automorphisms of complexes of curves and of Teichm\"uller spaces} : International Mathematics Research Notices (1997) 651-666.

[Kla] E. Klarreich, \tit{The boundary at infinity of the curve
complex and the relative Teichm\"uller space} : Preprint (1999).

[Lea] J. Leasure, \tit{Geodesics in the complex of curves of a surface} : PhD thesis (2002).

[MaMi1] H. A. Masur, Y. N. Minsky, \tit{Geometry of the complex of curves I: Hyperbolicity} : Inventiones Mathematicae \tbf{138} (1999) 103-149.

[MaMi2] H. A. Masur, Y. N. Minsky, \tit{Geometry of the complex of
curves II: Hierarchical structure} : Geometry \& Functional
Analysis \tbf{10} (2000) 902-974.

[PenH] R. C. Penner, J. L. Harer, \textit{Combinatorics of train tracks} : Annals of Mathematical Studies, Princeton University Press (1992).

[R] Mary Rees, \tit{The geometric model and large Lipschitz
  equivalence direct from Teichm\"uller geodesics} : Preprint (2004).
\end{document}